\documentclass{amsart}
\usepackage{amssymb}
\usepackage{xypic}

\newtheorem{theorem}{Theorem}[section]

\newtheorem{prop}[theorem]{Proposition}

\def\Hy{{\mathcal H}}

\def\cO{{\mathcal O}}

\def\N{{\mathbb N}}
\def\Q{{\mathbb Q}}
\def\R{{\mathbb R}}
\def\C{{\mathbb C}}
\def\F{{\mathbb F}}

\def\Aut{\operatorname{Aut}}
\def\Comm{\operatorname{Comm\,_H}}
\def\Isom{\operatorname{Isom}}
\def\Ker{\operatorname{Ker}}
\def\Stab{\operatorname{Stab}}

\def\POm{\operatorname{P\Omega}}
\def\PO{\operatorname{PO}}
\def\O{\operatorname{O}}
\def\SO{\operatorname{SO}}
\def\PSO{\operatorname{PSO}}
\def\SU{\operatorname{SU}}
\def\SL{\operatorname{SL}}
\def\PSL{\operatorname{PSL}}
\def\det{\operatorname{det}}

\def\LpF{{\mathcal L_p(F)}}
\def\M{{\mathfrak M}}
\def\NA{N_\Gamma(A)}

\begin{document}
\title{Finite Groups and Hyperbolic Manifolds}
\author{
Mikhail Belolipetsky
and
Alexander Lubotzky$^*$}
\thanks{$^*$ Partially supported by BSF (USA -- Israel) and ISF}
\address{M.~Belolipetsky \newline\phantom{iiii}
Institute of Mathematics,
Hebrew University,
Jerusalem 91904,
ISRAEL \newline\phantom{iiii}
Institute of Mathematics,
Koptyuga 4,
630090 Novosibirsk, RUSSIA
}
\email{mbel@math.nsc.ru}
\address{A.~Lubotzky \newline\phantom{iiii}
Institute of Mathematics,
Hebrew University,
Jerusalem 91904,
ISRAEL
}
\email{alexlub@math.huji.ac.il}

\begin{abstract}
The isometry group of a compact $n$-dimensional
hyperbolic manifold is known to be finite. We show that for every
$n\ge 2$, every finite group is realized as the full isometry
group of some compact hyperbolic $n$-manifold. The cases $n=2$ and
$n=3$ have been proven by Greenberg~\cite{G} and Kojima~\cite{K},
respectively. Our proof is non constructive: it uses counting results
from subgroup growth theory to show that such manifolds exist.

\end{abstract}

\maketitle

\section{Introduction}
\markboth{MIKHAIL BELOLIPETSKY AND ALEXANDER LUBOTZKY}{FINITE GROUPS AND HYPERBOLIC MANIFOLDS}

Let $\Hy^n$ denote the hyperbolic $n$-space, that is the unique connected,
simply connected Riemanian manifold of constant curvature $-1$. By a compact
hyperbolic $n$-manifold we mean a quotient space $M = \Gamma\backslash\Hy^n$
where $\Gamma$ is a cocompact torsion-free discrete subgroup of the group
$H = \Isom(\Hy^n)$ of the isometries of $\Hy^n$. The group $\Isom(M)$ of the
isometries of $M$ is finite and it is isomorphic to $N_H(\Gamma)/\Gamma$
where $N_H(\Gamma)$ denotes the normalizer of $\Gamma$ in $H$.

In 1972, Greenberg~\cite{G} showed that if $n=2$, then for every finite
group $G$ there exists a compact $2$-dimensional hyperbolic manifold $M$
(equivalently, cocompact $\Gamma$ in $H$) such that $\Isom(M)\cong G$
(equivalently, $N_H(\Gamma)/\Gamma \cong G$). A similar result for $n=3$
was proved in 1988 by Kojima~\cite{K}, who also mentioned the general
conjecture. The methods of Greenberg and Kojima are very much
of low dimensional geometry (Teichm\"uller theory and Thurston's Dehn
surgery, respectively).

The long standing problem of realizing every finite group as the isometry group
of some $n$-dimensional compact hyperbolic manifold is in flavor of the inverse
Galois problem and other questions of such kind (see e.g.~\cite{F}). What makes
the problem quite delicate is that even when it is solved for a group $G$,
it is still not settled either for the subgroups or for the factor
groups of $G$. In particular, our problem is non-trivial even for the
case of the trivial group, for which it means the existence of asymmetric
hyperbolic $n$-manifolds. Recently, Long and Reid~\cite{LR2}
showed that for every $n$ there exists a compact hyperbolic $n$-dimensional
manifold $M$ with $\Isom^+(M) = \{e\}$. Here, $\Isom^+(M)$ is the group of
orientation preserving isometries, it is a subgroup of index at most
two in $\Isom(M)$. They asked separately (\cite{LR2}, $\S4.3$ and
$\S4.4$) whether such an $M$ exists with $\Isom(M) = \{e\}$ as well
as for the general $G$. In this paper we give a complete solution to
the problem. It turns out, indeed, that the proof is somewhat different
for the cases $G = \{e\}$ and $G\neq\{e\}$.
\vskip 1ex

Our main result is the following
\begin{theorem}
For every $n\ge 2$ and every finite group $G$ there exist infinitely many
compact $n$-dimensional hyperbolic manifolds $M$ with $\Isom(M) \cong G$.
\end{theorem}

Let us describe the line of the proof. We start with the Gromov and
Piatetski-Shapiro construction~\cite{GPS} of a non-arithmetic lattice $\Gamma_0$
in $\Isom(\Hy^n)$. These lattices are obtained by interbreeding two
arithmetic lattices and the construction, in particular, implies that $\Gamma_0$
is represented as a non-trivial free product with amalgam. By
Margulis' theorem [Mr, Theorem~1, p.~2] the non-arithmeticity of the lattice
implies that its commensurator $\Gamma = \Comm(\Gamma_0)$ is a maximal
discrete subgroup of $\Isom(\Hy^n)$, so for every finite index
subgroup $B$ of $\Gamma$, $N_H(B) = N_\Gamma(B)$. It therefore
suffices to find such a torsion-free $B$ with $N_\Gamma(B)/B\cong G$.

We begin the search for $B$ inside $\Gamma_0$ which enables us to use its
amalgamated structure. To this end we modify the argument of~\cite{L2}
to show that $\Gamma$ has a suitable finite index torsion-free normal
subgroup $\Delta$ which is mapped onto a free group $F = F_r$ on $r\ge 2$
generators with kernel $M$. We then apply ideas and results from subgroup
growth theory [LS] to prove that $\Delta$ has a finite index
subgroup $A$ with $N_\Delta(A)/A$ isomorphic to $G$. A crucial
point here is that $F$ has at least $k!$ subgroups of index $k$ but at
most $k^{cr\log_2 k}$ of them are normal in $F$, for some absolute constant
$c$. (This result was proved in~\cite{L3} using the classification of the
finite simple groups, but the version we need is somewhat weaker and can
be proved without the classification. So the current paper is classification
free!) Another interesting group theoretic aspect is the use along the way
of a result from~\cite{L1} asserting that an automorphism of a free group
preserving every normal subgroup of a $p$-power index must be inner.
(The only known proof of this result relies on the theory of pro-$p$
groups.)

Now another problem has to be fixed. While $N_\Delta(A)/A \cong G$,
$N_\Gamma(A)$ can be (and in fact in many cases it is) larger than
$N_\Delta(A)$. To deal with this issue we modify $A$ by replacing it
by a somewhat larger subgroup $B$ of $\Gamma$ for which indeed
$N_\Gamma(B)/B \cong G$. Two delicate points have to be overcome
on the way: First is controlling the normalizer; what makes the whole
proof difficult is the fact that "normalizer is not continuous"; even
a small change from $A$ to $B$ can change the normalizer
dramatically. The second point is to keep $B$ torsion-free just as $A$.
This is achieved by keeping $B$ inside a suitable principal congruence
subgroup.

The paper is organized as follows: In \S2 we collect a number of group
theoretic results to be used in the later sections. In \S3 we  bring the
main group theoretical method to find finite index subgroups $B$ in a
group $\Gamma$ with $N_\Gamma(B)/B \cong G$. Then in \S4, a group $\Gamma$
in $\Isom(\Hy^n)$ is constructed which satisfies all the needed assumptions.
So only \S4 contains some geometry. We end with remarks and suggestions
for further research in \S5.

\bigskip

{\bf\noindent Acknowledgement.} The first author is grateful to the Hebrew
University for its hospitality and support.
The second author thanks Alan Reid for sending the manuscript \cite{LR2}
which was the main inspiration for this work, and B.~Farb and S.~Weinberger
for a helpful discussion.

\section{Group theoretic preliminaries}

In this section we will present a number of group theoretic results which
we will use later. We begin with the {\it free groups}.

\subsection{}
Let $F = F_r$ be a free group on $r\ge 2$ generators. For a prime $p$ denote
by $\LpF$ the family of all normal subgroups of $F$ of $p$-power
index.

\begin{theorem}\cite{L1}
(a) $F$ is residually-$p$, i.e., $\bigcap_{N\in\LpF} N = \{e\}$.

\noindent
(b) If $\alpha$ is an automorphism of $F$ such that $\alpha(N) = N$ for
any $N\in\LpF$ then $\alpha$ is inner.
\end{theorem}
While part~(a) is well known and easy, it is interesting to remark that the
proof of part~(b) in~\cite{L1} is based on the work of Jarden~-~Ritter~\cite{JR}
which combines pro-$p$ groups and relation modules.

\begin{prop} {\rm (Schreier's theorem, cf. [S, Theorem~5, p.~29])}
If $H$ is a subgroup of $F = F_r$ of index $k$ then $H$ is a free group on
$1 + k(r-1)$ generators.
\end{prop}

For a finitely generated group $\Gamma$ we denote by $a_n(\Gamma)$ (resp.,
$a^\lhd_n(\Gamma)$, $a^{\lhd\lhd}_n(\Gamma)$) the number of subgroups (resp.,
normal subgroups, subnormal subgroups) of index $n$ in $\Gamma$, and let
$s^\bullet_n(\Gamma) = \sum_{i=1}^n a^\bullet_i(\Gamma).$


\begin{theorem}
(a) {\rm (cf. [LS, Corollary~1.1.2, p.~13 and Corollary~2.1.2, p.~41])}
$$(n!)^{r-1} \le a_n(F) \le n(n!)^{r-1}.$$

\noindent
(b) {\rm (cf. \cite{L3}, see also [LS, Theorem~2.6])} There exists a constant $c$
such that $$a_n^\lhd(F) \le n^{cr\log_2 n}.$$

\noindent
(c) {\rm (cf. [LS, Theorem~2.3])} $s_n^{\lhd\lhd}(F) \le 2^{rn}$.
\end{theorem}

\subsection{}
We now turn to {\it free products with amalgam} and {\it HNN-constructions}.
Let $Q$ be a finite group of order $\ge 3$ and $T$ be a subgroup of $Q$
satisfying the following:
\medskip

$(*)$ If $N\le T\le Q$ and $N\lhd Q$ then $N = \{e\}$.
\medskip

In particular, $(*)$ implies that $[Q:T]>2$.
\medskip

Let $R = Q*_T Q$ be the free product of $Q$ with itself amalgamated along $T$ or
$R = Q*_T$~--~the HNN-construction.
There is a natural projection $\bar\pi : R\to Q$ whose kernel will be denoted by
$F = \Ker(\bar\pi$). As $Q$ is also a subgroup of $R$ we have $R = F\rtimes Q$~--~a
semi-direct product.

\begin{prop}
(a) $F$ is a non-abelian free group.

\noindent
(b) $C_R(F) = \{e\}$, i.e., the centralizer of $F$ in $R$ is trivial.

\noindent
(c) If $\alpha$ is an automorphisms of $R$ satisfying $\alpha(F) = F$ and
$\alpha |_F = id$ then $\alpha = id$.
\end{prop}

\begin{proof}
(a) 
This is a known fact, let us briefly recall the argument. By definition $F\smallsetminus\{e\}$
does not meet any conjugate of $Q$, so it acts trivially on the tree associated
to $R$ by the Bass-Serre theory [S, Section 4.2]. This implies that $F$ is
a free group [S, Theorem~4, p.~27]. The rank $r$ of $F$ can be computed using
the formula from [S, Exercise~3, p.~123], in particular, the condition $[Q:T] > 2$
implies that $r\ge 2$ and so $F$ is a non-abelian free group.

(b) The group $C=C_R(F)$ is a normal subgroup of $R$.
As $F$ is a non-abelian free group with a trivial center,
$C\cap F = \{e\}$ and hence $C$ is finite.
We claim that $C\le T$. If $R$ is a free product with amalgamation, then by [S, Theorem~7, p.~32],
$R$ acts on a tree with a fundamental domain consisting of an edge $e$ with two vertices $v_1$
and $v_2$ such that $\Stab_R(v_i)$, $i = 1,2$, are two copies of $Q$ in $R$
which we denote by $Q_1$ and $Q_2$. Since $C$ is finite it is conjugate into one
of the $Q_i$'s, say $Q_1$. But since $C$ is normal in $R$ it is contained in $Q_1$
and in all the conjugates of $Q_1$ in $R$. This implies that $C$ fixes all the
vertices in the $R$-orbit of $v_1$. The latter fact also implies that $C$ fixes
$v_2$ (and all the vertices in its orbit). Hence $C$ is also in $Q_2$ and so
$C\le Q_1\cap Q_2 = T$. The case of HNN-construction is even easier. Again $R$
acts on a tree and this time the fundamental domain contains a single vertex.
We get that $C$ is contained in the stabilizers of all the vertices and edges,
so $C$ is contained in $T$. In both cases condition $(*)$ implies that $C=\{e\}$.

(c) Let $f\in F$ be an arbitrary element of $F$ and $q\in Q$ an arbitrary element
of $Q$. We have $q^{-1}fq\in F$, so $\alpha(q^{-1}fq) = q^{-1}fq$ and $\alpha(f)=f$.
This implies that
$$ q^{-1}fq = \alpha(q^{-1}fq) = \alpha(q)^{-1}f\alpha(q), $$
and hence
$$ \alpha(q)^{-1}q \in C_R(F) = \{e\}.$$
Thus $\alpha(q) = q$ and so $\alpha$ is the identity on $Q$. Since $R = F\rtimes Q$,
$\alpha$ is the identity automorphism of $R$.
\end{proof}

\subsection{}
Finally, we will need some facts on the {\it finite orthogonal groups} (see e.g. [A]).
Let $f$ be an $m$-dimensional quadratic form over a finite field $\F$ of characteristic
$p>2$. If $m$ is odd, there is a unique, up to isomorphism, orthogonal group $\O(f)=\O_m(\F)
= \O_m$. If $m$ is even there are two groups $\O^+_m$ and $\O^-_m$ corresponding to
the cases when $f$ splits and does not split over $\F$, respectively. Let $\SO(f)$,
$\PSO(f)$ denote the corresponding special orthogonal and projective special orthogonal
groups, and let $\Omega(f) = [\O(f),\O(f)]$ be the commutator subgroup of $\O(f)$.
The projective group $\POm(f) = \POm^\pm_m(\F)$ is generally simple and is contained in $\PSO(f)$
with index at most $2$. More precisely, $\POm_m$ is simple if $m\ge5$ or $m=3$ and $p>3$;
in case $m=4$, $\POm^-_4$ is simple but $\POm^+_4 = \POm_3\times\POm_3$ is a direct
product of two groups which are simple if $p>3$; the cases $m=1,2$ will not be used in
this paper. We sometimes omit the $\pm$ sign is the notations.
For future reference note also that the centralizer of $\POm(f)$ in $\PO(f)$
is always trivial.

\section{The main algebraic result}

In this section we prove a purely group theoretic result. In Section~4 we will
show that it can be implemented for suitable non-arithmetic lattices in $\PO(n,1)$.

Throughout this section $\Gamma$ is a finitely generated group, $\Delta$ a finite
index normal subgroup of $\Gamma$, and $M$ is a normal subgroup of $\Delta$ with
$\Delta/M$ being isomorphic to a free group $F=F_r$ on $r\ge 2$ generators. We
denote by $N = N_\Gamma(M)$ the normalizer of $M$ in $\Gamma$, so
$\Delta\le N\le\Gamma$. The group $N$ acts by conjugation on $F = \Delta/M$. Denote
by $C = C_N(\Delta/M)$ the kernel of this action and by $D$ the subgroup of all
elements of $N$ which induce inner automorphisms of $F$. Both $C$ and $D$ are
normal in $N$ and clearly
$$D = \Delta C.$$
Moreover, $M$ is normal in $D$, $\Delta$ and in $C$. As $F = \Delta/M$ has a
trivial center and hence intersects $C/M$ trivially, $\Delta\cap C = M$. So,
taking mod $M$ we get
$$D/M = \Delta/M \times C/M.$$
Moreover, $\Delta/M$ is of finite index in $D/M$ and hence $C/M$ is a finite group.
$$
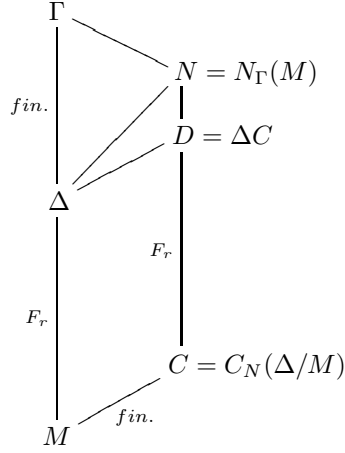
\xymatrix@C=2pt@R=2.8ex{
*+[r]{\Gamma} \ar@{-}[ddd]_{fin.} \ar@{-}[dr] & \\
                                &*[r]{\:N = N_\Gamma(M)} \ar@{-}[d] \ar@{-}[ddl]  \\
                                &*+[r]{D = \Delta C} \ar@{-}[dl] \ar@{-}[dddd]_{F_r}\\
*+[r]{\Delta} \ar@{-}[dddd]_{F_r}    & \\
                                & \\
                                & \\
                                &*+[r]{C = C_N(\Delta/M)} \ar@{-}[ld]^{fin.} \\
{M} \\ }
$$



\vskip 1ex

This section is devoted to the proof of the following result.

\begin{theorem}
Let $\Gamma$, $\Delta$ and $D$ be as above. For every finite group $G$ there exist
infinitely many finite index subgroups $B$ of $D$ with $N_\Gamma(B)/B$ isomorphic
to $G$.
\end{theorem}

\begin{proof}
Denote the order of $G$ by $g = |G|$, and $g' = |\Aut(G)|$. Also let $d+1 = [\Gamma :D]$
and $e,\gamma_1,\dots,\gamma_d$ be representatives of the right cosets of $D$ in
$\Gamma$, i.e. $\Gamma\smallsetminus D = \bigcup_{i=1}^d D\gamma_i$.

Let $x > max\{g,\: g'\}$ be a very large integer, to be determined later. Choose $d+1$
primes $p_0 < p_1 < \ldots < p_d$ with $p_0 \ge x^2$.

Now, if $d = 0$, choose a normal subgroup of $\Delta$ of index $p_0$ containing $M$ and
call it $K$. If $d>0$ the definition of $K$ will be more delicate: We claim that for every
$i = 1,\ldots,d$ there exists a normal subgroup $K_i\supset M$ of index $p_i^{\alpha_i}$ in
$\Delta$ for some $\alpha_i\in\N$, such that $K_i^{\gamma_i} \neq K_i$ (where
$K_i^{\gamma_i} = \gamma_i^{-1}K_i\gamma_i$). Indeed, if not then conjugation by
$\gamma_i$ stabilizes all normal subgroups of $\Delta$ containing $M$ and of index
$p_i$-power in $\Delta$. As $\Delta/M\cong F$ is residually-$p_i$ (Theorem~2.1(a)),
$M$ is the intersection of these normal subgroups and hence $\gamma_i$ normalizes
$M$ and so $\gamma_i\in N_\Gamma(M)$. Moreover, $\gamma_i$ acting on $F = \Delta/M$
is now an automorphism of $F$ which preserves any normal subgroup of $F$ of
{\nobreak $p_i$-power} index. By Theorem~2.1(b), $\gamma_i$ induces on $F$ an inner
automorphism and hence $\gamma_i\in D$ in contradiction to the way $\gamma_i$
was chosen. We define $K := \bigcap_{i=1}^d K_i$.
\medskip

In both cases denote the index of $K$ in $\Delta$ by $k$. Observe that $k\ge x^2$
and $K/M$ is a subgroup of $\Delta/M = F = F_r$ of index $k$ so by Proposition~2.2,
$K/M$ is a free group on $1+k(r-1)$ generators.
\medskip

The proof now splits into two cases.
\medskip

{\noindent\bf Case 1:} $G\neq\{e\}$. We claim that there are at least
$g^{1+k(r-1)-\log_2 g}$ different epimorphisms from $K$ onto $G$ whose kernels contain
$M$. Indeed, it follows from an easy argument that a finite group of order $g$ is generated
by at most $\log_2 g$ elements (cf.~[LS, Lemma~1.2.2, p.~14]), so we can send the first $\log_2 g$
generators of $K/M \cong F_{1+k(r-1)}$ to fixed generators of $G$ and all the rest generators
of $K/M$ can be sent arbitrarily into $G$. Two such epimorphisms $\phi_1$ and $\phi_2$
have the same kernel if and only if there exists $\beta\in\Aut(G)$ such that
$\phi_1 = \beta\circ\phi_2$. Again, as $G$ is generated by $\log_2g$ elements,
$g' = |\Aut(G)| \le g^{\log_2g}$. Hence there exist at least
$$\frac{1}{g'}\:g^{1+k(r-1)-\log_2g} \ge g^{k(r-1)-2\log_2g} =: z$$
normal subgroups $A$ of $K$ containing $M$ with $K/A \cong G$.

Let $\M$ denote the set of these subgroups $A$.
We claim that for every $A\in\M$, $\NA\le D$.
If $d=0$, $\Gamma=D$ and there is nothing to prove. If $d>0$, let $\gamma\in\NA$
and assume $\gamma\in\Gamma\smallsetminus D$, so $\gamma = \delta\gamma_i$ for some
$\delta\in D$ and some $i\in\{1,\dots,d\}$. Now, $\gamma$ normalizes $\Delta$
and $A$, hence $\gamma$ normalizes $K_i$ since $K_i$ is the only normal subgroup
of $\Delta$ of index $p_i^{\alpha_i}$ containing $A$. Moreover, $\delta$ being
an element of $D$ induces an inner automorphism on $F=\Delta/M$ hence it also
normalizes $K_i$. Thus $\gamma_i = \delta^{-1}\gamma$ normalizes $K_i$, but this
contradicts the way $K_i$ was chosen. So for every $A\in\M$, $\NA\le D$.

Let us observe that $KC\le\NA$ (since $K$ and $C$ are normal subgroups of $D$,
their product $KC$ is indeed a subgroup). We will show next that for {\it some}
$A\in\M$, $\NA = KC$. This will be done by a counting argument.

We already know that $\NA = N_D(A) \supseteq KC$. As $D/M \cong \Delta/M \times C/M$,
$N_D(A)/M \cong N_\Delta(A)/M \times C/M$, so we
can project everything to $\Delta/M$ and it suffices to show that there exists $A\in\M$
with $N_\Delta(A) = K$.

Fix one $A\in\M$ and denote $L := N_\Delta(A)$, $l := [\Delta : L]$.
The subgroup $L$ contains $M$ and $L/M$ is a subgroup of $\Delta/M \cong F_r$ of index $l$,
so $L/M$ is a free group on $1+l(r-1)$ generators. Also, $A$ is a normal subgroup of $L$ of
index $gk/l$, so $A/M$ is a normal subgroup of $L/M$ of the same index. By Theorem~2.3(b),
$L/M$ has at most
$$(\frac{gk}{l})^{c(1+l(r-1))\log_2(gk/l)} \le 2^{crl\log_2^2(gk/l)} =: y$$
normal subgroups of index $gk/l$, where $c$ is an absolute constant.
Since $L\ge K$, $l$ divides $k$ and it is a proper divisor
if $L\neq K$. So $L\gneqq K$ implies that $l$ is at most $k/p_o \le k/x^2 < k/x$.
Note also that as $k\ge x^2$, $k/x \ge x$.
Running over all $l$ in the range between $1$ and $k/x$, for big enough $x$,
the maximal value of $y$ is obtained when $l = k/x$. Thus, there are at most
$2^{cr(k/x)\log_2^2(gx)}$ subgroups in $\M$ which are normalized by a given $L \ne K$.

Now note that $\Delta/K$ is a finite nilpotent group
(of order $p_0$ if $d=0$ and of order $\Pi_{i=1}^d p_i^{\alpha_i}$ otherwise), $L/K$
is a subgroup of $\Delta/K$ and any subgroup of a nilpotent group is subnormal. This
implies that $L$ is a subnormal subgroup of $\Delta$ of index less than $k/x$ (if $L\neq K$),
hence by Theorem~2.3(c) there are at most $2^{rk/x}$ possibilities for such $L$'s.

Putting all this together we see that there are at least $z = g^{k(r-1)-2\log_2 g}$
possibilities for $A$ and out of them at most
$$w := 2^{rk/x}2^{cr(k/x)\log_2^2(gx)} = 2^{(rk/x)(1+c\log_2^2(gx))}$$
have their normalizer bigger then $KC$. For $x$ large enough,
$z > w$ and so there exist $A\in\M$ with $\NA = KC$. In fact, most $A\in\M$ do
satisfy this.

Note however, that for these $A$
$$\NA/A = KC/A \cong K/A \times C/M \cong G \times C/M.$$
So, we have not achieved the goal of Theorem~3.1 yet. For this we will have to enlarge
$A$, but first let us consider the second case.
\bigskip

{\noindent\bf Case 2:} $G =\{e\}$. Let $q$ be a prime close to $k$ and different from
$p_0,\ldots,p_d$, such that $q > p_0 \ge x^2$.
Let $\M$ be the set of all subgroups of $K$ of index $q$ containing $M$. By Theorem~2.3(a)
there are at least
$$(q!)^{k(r-1)+1} > (q!)^{k(r-1)} =: z $$
such subgroups. As before we claim that for every $A\in\M$, $\NA\le D$. If $d =0$ there is
nothing to prove. If $d > 0$, the argument is exactly the same: write
$\gamma\in\Gamma\smallsetminus D$ as $\gamma = \delta\gamma_i$ for some $\delta\in D$ and
$i\in\{1,\ldots,d\}$, if $\gamma$ normalizes $A$ then it also normalizes $K_i$, which implies
that $\gamma_i$ normalizes $K_i$~-- a contradiction.

So, again for every $A\in\M$, $\NA\le D$. This time we want to prove that for some (in fact,
for most) $A\in\M$, $\NA = AC$ (note: not $KC$ as for $G\neq\{e\}$).
As in the previous case we
will project everything to $\Delta/M$ and show that there exists $A\in\M$
with $N_\Delta(A) = A$.

To this end, note that if $L:=N_\Delta(A)$ is strictly larger than $A$, then $L$ is
a subgroup of index, say, $l$ in $\Delta$, $l$ being a proper divisor of $kq$ and so
$l\le kq/x$. Similarly to the previous case $L/M$ is a free group on $1+l(r-1) < lr$
generators and $A/M$ is its normal subgroup of index $kq/l$. It follows from Theorem~2.3(b)
that $L$ can normalize at most
$$(\frac{kq}{l})^{clr\log_2(kq/l)} = 2^{crl\log_2^2(kq/l)} =: y$$
subgroups from $\M$. Now, $l$ can take on values between $1$ and $kq/x$, and clearly the
maximum of $y$ is attained when $l = kq/x$ (if $x$ is big enough) corresponding to
$y = 2^{cr (kq/x)\log_2^2x}$.

Continuing to work mod $M$, $L/M$ is a subgroup of index at most $kq/x$ of the free group
$\Delta/M \cong F_r$. There are, therefore, at most $((kq/x)!)^r$ 
such subgroups by Theorem~2.3(a). The latter number is bounded by
$(kq/x)^{c'(kq/x)r}$ 
for a suitable constant $c'$. Altogether, at most
$$ w := (kq/x)^{c'(kq/x)r}\: 2^{cr(kq/x)\log_2^2x}
= 2^{r(kq/x)(c'\log_2(kq/x) + c\log_2^2x)} $$
of $A \in \M$ have normalizer $N_\Delta(A)$ larger then $A$. Recall that $q$
was chosen to be approximately $k$ and $z = (q!)^{(r-1)k}$. An easy estimate shows
that $z > w$ provided $x$ is large enough. Thus, there exists $A$ in $\M$ (in fact,
most $A\in\M$) for which $N_\Delta(A) = A$ and $\NA = N_D(A) = AC$.

\bigskip

Let us now treat both cases $G \neq \{e\}$ and $G = \{e\}$ together: The above
arguments show that we can always find a subgroup $A$ of a finite index in
$\Delta$ such that $N_\Delta(A)/A \cong G$, $\NA = N_D(A)$ and
$\NA/A \cong G\times C/M$. Let us replace $A$ by $B = AC$. Since $A$ and $C$
are both in $D$ and $C\lhd D$, $B$ is indeed a subgroup and it is contained
in $D$. It is also clear that $B\cap\Delta = A$ (look at everything mod $M$!).
We claim that $N_\Gamma(B)/B \cong G$. This will finish the proof of the
theorem.

First, note that if
$\gamma\in N_\Gamma(B)$ then $\gamma$ also normalizes $\Delta$ (since $\Delta\lhd\Gamma$)
and hence:
$$A^\gamma = (B\cap\Delta)^\gamma = B^\gamma\cap\Delta^\gamma = B\cap\Delta = A,$$
so $\gamma\in\NA$. On the other hand every $\gamma\in\NA$ also normalizes $C$,
since $\NA\le D$ and $C$ is normal in $D$. This shows that $N_\Gamma(B) =\NA$
and so
$$N_\Gamma(B)/B = \NA/AC \cong N_\Delta(A)/A \cong G$$
as claimed.

\bigskip

We finally mention that by choosing infinitely many different $x$'s (and hence also
the $p_i$'s) we will get infinitely many subgroups $B$ of $\Gamma$
with $N_\Gamma(B)/B \cong G$.
\end{proof}

\section{Geometric realization}

In this section we will show that for every $n\ge 2$ there exist
non-arithmetic lattices in $H = \Isom(\Hy^n)$ satisfying the
assumptions of Theorem~3.1 and then deduce the main result of
this paper. Recall that $H$ can be identified with $\O_0(n,1)$~--
the subgroup of the orthogonal group $\O(n,1)$ which preserves
the upper-half space, it is isomorphic to the projective
orthogonal group $\PO(n,1) = \O(n,1)/\{+1,-1\}$. The subgroup
$\SO_0(n,1)$ of $\O_0(n,1)$ of all the elements of $H$ with
determinant $1$, is the group of orientation preserving isometries.

\begin{prop} For every $n\ge 2$ there exist a maximal cocompact
non-arithmetic lattice $\Gamma$ in $H$ with subgroups $M$, $\Delta$
and $D$ satisfying the following:
\begin{itemize}
\item[(i)] $\Delta\lhd\Gamma$ and $[\Gamma : \Delta] < \infty$.
\item[(ii)] $M\lhd\Delta$, $\Delta/M$ is a non-abelian free group.
\item[(iii)] $[\Gamma : D] < \infty$, $D\le N_\Gamma(M)$ and
\newline $D = \{\delta\in N_\Gamma(M) \mid \delta {\text\ induces\ an\
inner\ automorphism\ on\ } \Delta/M \}.$
\item[(iv)] $D$ is torsion-free.
\end{itemize}
\end{prop}

\begin{proof}
Let $\Gamma_0$ be a cocompact non-arithmetic lattice in $H$
obtained by Gromov~-~Pia\-tet\-ski-Shapiro construction~\cite{GPS}.
Recall that $\Gamma_0$ is constructed as follows:
One starts with two non-commensurable torsion-free arithmetic lattices
$L_1$ and $L_2$ in $H$,
such that each of the corresponding factor manifolds
$W_i = L_i\backslash\Hy^n$ admits a totally geodesic hypersurface
$Z_i$ ($i = 1,2$) and $Z_1$ is isometric to $Z_2$. Assume that
$Z_i$ ($i=1,2$) separates $W_i$ into two pieces $X_i\cup Y_i$ (the
non-separating case can be treated in a similar way). Then a new
manifold $W$ is defined by gluing $X_1$ with $Y_2$ along $Z_1$
(which is isomorphic to $Z_2$). In particular, $W$ itself has a properly
embedded totally geodesic hypersurface $Z$ (isometric to $Z_1$ and $Z_2$)
and so $\pi_1(W) = \pi_1(X_1)*_{\pi_1(Z)} \pi_1(Y_2)$. As explained
in~\cite{GPS}, $\Gamma_0 = \pi_1(W)$ is a non-arithmetic lattice in
$H  \cong \O_0(n,1)$ which can be supposed to be contained in
$\SO_0(n,1)$ (so $W$ is orientable) and $\pi_1(Z)$ is a subgroup (in fact,
a lattice) in a conjugate of $\SO_0(n-1,1)$.

Let $\Gamma = \Comm(\Gamma_0) = \{ g\in H \mid
[\Gamma_0 : \Gamma_0 \cap g^{-1}\Gamma_0g] < \infty\}$ be the
commensurability group of $\Gamma_0$. Since $\Gamma_0$ is
non-arithmetic, Margulis Theorem~[Mr, Theorem~1, p.~2] implies that
$\Gamma$ is also a lattice, a maximal lattice in $\Isom(\Hy^n)$.

Let now $\Lambda$ be a finite index normal subgroup of $\Gamma$
which is contained in $\Gamma_0$. So $\Lambda$ is the fundamental
group of a finite sheeted cover $W'$ of $W$. Pulling back the
hypersurface $Z$ to $W'$, we deduce that $W'$ also admits a properly
embedded totally geodesic hypersurface and hence $\Lambda
= \Lambda_1*_{\Lambda_3}\Lambda_2$ or $\Lambda = \Lambda_1*_{\Lambda_3}$
is a non-trivial free product with amalgam or an HNN-construction.

By [GPS, Corollary~1.7.B] the groups $\Lambda_1$ and $\Lambda_2$ are Zariski dense
in $\SO(n,1)$. By the construction these groups are contained in $\Lambda\le\Gamma$.
Let $\cO$ (resp., $\cO_i$, for $i = 1, 2$) be the minimal ring of definition of $\Gamma$
(resp., $\Lambda_i$) in the sense of~\cite{V}. As $\Lambda$ is finitely generated
group, $\cO$ is a finitely generated ring. In fact,
it is contained in some number field $k$. This last claim is true
for all lattices in $\O(n,1)$ if $n\ge 3$ by the local rigidity of these lattices
(see [R, Proposition~6.6, p.~90]). But it also follows (for every $n$, including
$n=2$) for the Gromov~-~Pia\-tet\-ski-Shapiro lattices directly from their
construction.

Thus $\cO$ is a ring of $S$-integers in some (real) number field $k$. For
$i=1, 2$, $\cO_i$ is a subring of $\cO$, so it is the ring of $S_i$-integers
of some subfield $k_i$ of $k$ for a suitable finite set $S_i$ of primes
in $k_i$. We can assume $\Gamma\le\SO(n,1)(\cO)$ and $\Lambda_i\le\SO(n,1)(\cO_i)$
for $i=1, 2$.

Now, the strong approximation for linear groups ([We, Theorem~1.1], [P, Theorem~0.2],
see also [LS, Window~9])
implies that for almost every maximal ideal $\mathcal P$ of $\cO$ with
finite quotient field $\F_q = \cO/\mathcal P$, $q = |\cO/\mathcal P|$, the image
of $\Gamma$ in $\PO_{n+1}(\F_q)$ contains $\POm_{n+1}(\F_q)$ which is of index
at most two in $\PSO_{n+1}(\F_q)$.
The same also applies to $\Lambda$ since
$\Lambda$ is of finite index in $\Gamma$. Moreover $\Lambda_i$ is also Zariski
dense for $i=1, 2$, so a similar statement holds for $\Lambda_i$ with respect to the
ring $\cO_i$. By Chebotarev density theorem, there exist infinitely many primes $l$
in $\Q$ which split completely in $k$ (and hence also in $k_i$). Thus for every
prime ideal $\mathcal P$ of $\cO$ which lies above such $l$,
$\cO/\mathcal P = \cO_i/\cO_i\cap\mathcal P \cong \F_l$. Moreover, if we replace
$\Lambda$ by the intersection of all the index $2$ subgroups in it (note
that this intersection is characteristic in $\Lambda$ and so normal in
$\Gamma$), we can assume that for infinitely many rational primes $l$, the
images of $\Lambda$, $\Lambda_1$ and $\Lambda_2$ are exactly the groups
$\POm_{n+1}(\F_l)$.

Choose such a prime $l$. We obtain a homomorphism
$$\pi:\Lambda\to Q = \POm_{n+1}(\F_l)$$
with $\pi(\Lambda_1) =\pi(\Lambda_2) = Q$ while $T = \pi(\Lambda_3)\le \PSO_n(\F_l)$
is a proper subgroup of $Q$ (and by choosing $l$ sufficiently large we
can assume that the index of $T$ in $Q$ is as large as we want).
For later use we observe that if $T$ contains a normal subgroup $N$ of $Q$ then $N=\{e\}$.
Indeed, if $n\neq 3$ or $n=3$ and $Q = \POm^-_4$, $Q$ is a finite simple group and there
is nothing to prove. Suppose $n=3$, $Q = \POm^+_4 \cong \POm_3\times\POm_3$. The only
possibility for $N\neq\{e\}$ is $N = \POm_3$. The image $T$ of $\Lambda_3$ in $\Omega_4$
is equal to $\Stab_{\Omega_4}(U)$,
the stabilizer of a $3$-dimensional subspace $U$ of $V = \F_l^4$ and hence indeed it is
isomorphic to $\Omega_3$, but it cannot be a normal subgroup of $\Omega_4$.
For if this is the case, then for every $g \in \Omega_4$,
$\Stab_{\Omega_4}(gU) = gTg^{-1} = T$. This implies $gU = U$ (otherwise $\Omega_3$ would
preserve a $2$-dimensional subspace). Now, this means that $U$ is
$\Omega_4$ invariant, which is a contradiction. (We recall that
$\POm^+_4 \cong \POm_3\times\POm_3 \cong \PSL(2)\times\PSL(2)$ via the
action of $\SL(2)\times\SL(2)$ on the $2\times 2$ matrices by $(g,h)(A) = gAh^{-1}$,
with $\det(A)$ being the invariant quadratic form. But, the action on the
$4$-dimensional space is irreducible.)

The universal property of free products with amalgam and HNN-constructions
implies that there exist a homomorphism
$$ \tilde\pi: \Lambda \to R  = Q*_TQ\ ({\rm or}\ =Q*_T)$$
depending on $\Lambda = \Lambda_1*_{\Lambda_3}\Lambda_2$ or $\Lambda = \Lambda_1*_{\Lambda_3}$.
The group $R$ is mapped by $\bar\pi$ onto $Q$ with a kernel $F$
which is a non-abelian free group by Proposition~2.4(a). Let $M = \Ker(\tilde\pi)$
and $\Delta = \Ker(\bar\pi\circ\tilde\pi)$. We have:
$$ 
\xymatrix@R=2.8ex{
& & R \ar[dr]_{\bar\pi}\\
\Delta \ar[r] & \Lambda \ar[ur]_{\widetilde\pi} \ar[rr]_{\pi} & & Q \\
M \ar@{-}[u] \ar[ur] \\ }
$$
It is easy to see that $M\lhd\Delta$
and $\Delta/M \cong F$. Also, $\Delta = \Lambda\cap\Gamma(l)$ where
the congruence subgroup $\Gamma(l)$ is the kernel of the projection of
$\Gamma$ to $\PO_{n+1}(\F_l)$. Thus $\Delta$ is  a finite index normal
subgroup of $\Gamma$. We therefore have the properties~(i) and (ii) of the
proposition.

Let now $D = \{\delta\in N_\Gamma(M) \mid \delta {\rm\ induces\ an\
inner\ automorphism\ on\ } F\cong\Delta/M \}.$ Since
$D$ contains $\Delta$, we are left only with proving that $D$ is
torsion-free. Denote $C = \{ \delta\in N_\Gamma(M) \mid
\delta|_{\Delta/M} = id \}$. Then $D = \Delta C$. We will show
that $C\le\Gamma(l)$ which will prove that $D\le\Gamma(l)$. As
$\Gamma(l)$ is torsion-free (when $l$ is large enough), we will
deduce that $D$ has no torsion.

Let $c\in C$. The element $c$ acts on $\Lambda/M\cong R$ with
restriction to $\Delta/M\cong F$ being trivial. Such an
automorphism of $R$ is trivial by Proposition~2.4(c), so
$c$ acts trivially on $\Lambda/M$, i.e. $[c,\Lambda]\subseteq M$.
Taking this mod $\Gamma(l)$ we deduce that $c$ centralizes
$Q = \POm_{n+1}(\F_l)$ in $\PO_{n+1}(\F_l)$, so $c$ is trivial there by
Section~2.3. This implies that $c\in\Gamma(l)$ and Proposition~4.1 is
now proved.
\end{proof}

We can now prove the main Theorem from the Introduction:
\bigskip

Let $n\ge2$, $G$ is a finite group, $\Gamma$ and $D$ are lattices
as in Proposition~4.1. By Theorem~3.1, there exist infinitely many
finite index subgroups $B$ of $\Delta$ with $N_\Gamma(B)/B\cong G$.
Note that $N_H(B) \le \Comm(\Gamma) = \Gamma$, so $N_H(B) = N_\Gamma(B)$
and hence $N_H(B)/B\cong G$. As explained in the introduction
$N_H(B)/B \cong \Isom(B\backslash\Hy^n)$. The theorem is proved.

\section{Remarks}

\subsection{}
As already pointed out in the Introduction, although in Theorem~2.3(b)
we make use of the results of~\cite{L3}, which rely on the classification
of the finite simple groups, what we really need in this paper does not require
the classification. Indeed, Theorem~2.3(b) says that
$a_n^\lhd(F_r) \le n^{cr\log_2 n}$, where $a_n^\lhd(F_r)$ denotes the number
of index $n$ normal subgroups $N$ of the free group $F = F_r$ and $c$ is a
constant. But when we use it in the proof of Theorem~3.1, we need such an
upper bound only for those $N\lhd F$, for which $F/N$ has a normal subgroup
$G_0$ isomorphic to the fixed finite group $G$ and $(F/N)/G_0$ is nilpotent.
In~\cite{Mn}, Mann showed implicitly that if $\mathcal S$ is a family of finite
simple groups such that for every $S\in\mathcal S$, $S$ has a presentation
with at most $c_0(\mathcal S)\log_2|S|$ relations, then there is a constant
$c_1(\mathcal S)$ such that the number of index $n$ normal subgroups $N$ of
$F_r$ with all composition factors of $F/N$ being from $\mathcal S$ is
bounded by $n^{c_1(\mathcal S)r\log_2 n}$. Mann's argument is elementary.
We could use this result instead of Theorem~2.3(b). Since all the composition
factors of our groups are either those of $G$ or abelian, clearly, such a
$c_0(\mathcal S)$ exists.

\subsection{}
The finite volume non-compact case can be treated in an entirely similar way.
So, for every finite group $G$ there also exist infinitely many finite volume
non-compact $n$-dimensional hyperbolic manifolds $M$ with $\Isom(M)\cong G$.

\subsection{}
The proof of Proposition~4.1 actually shows that the subgroup $D$ constructed
there is contained in a principle congruence subgroup $\Gamma(l)$. Now, if
$l>2$ (which was indeed an assumption), then $\Gamma(l)\subset\SO(n,1)$. So
we actually provided infinitely many $M$'s with $\Isom(M)=\Isom^+(M)\cong G$.

P.~M.~Neumann suggested the following generalization for the problem:
Let $G$ be a finite group with a subgroup $G^+$ of index $2$. For every $n \ge 2$
does there exist a compact hyperbolic $n$-manifold $M$ with an isomorphism
$\psi : \Isom(M)  \to G$ such that
$\psi(\Isom^+(M)) = G^+$ ?



\subsection{} Our argument is close in spirit to Greenberg's proof for $n=2$:
while he counts the  dimensions of certain subspaces of the Moduli spaces
we use counting results on subgroups growth which also allow to detect the
existence of the manifolds with the prescribed groups of symmetries.
The method of Long and Reid is constructive in a sense.

Our method is not constructive, still the proof says something about its effectiveness.
For a fixed $n$, one has to find $\Gamma_0$~-- the Gromov~-~Piatetski-Shapiro
lattice in $\O(n,1)$ as in \S4. Then, with the notations of \S4, we need to find
a prime $l$ for which the image of $\Lambda$
(and also of its subgroups $\Lambda_1$ and $\Lambda_2$) to $\PO_{n+1}(\F_l)$
contains $\POm_{n+1}(\F_l)$ (all but finitely many primes which split in the
ring of definition of $\Lambda$ have this property). Once this is done the
proof gives an explicit estimate for the index of $B$ in $\Lambda$ for which
$N_H(B)/B\cong G$.

One can define for a finite group $G$, $f(n,G)$ to be the minimal volume of a
compact $n$-dimensional manifold $M$ with $\Isom(M) \cong G$. It may be of interest
to give some bounds on $f(n,G)$.

We mention by passing that for $n\ge 4$ and a given $r > 0$, there are only
finitely many $n$-dimensional hyperbolic manifolds of volume at most $r$~\cite{Wa}.
In~\cite{BGLM}, it is shown that the growth rate of the number of manifolds
is like $\exp(c(n)r\log r)$. One may ask, whether for {\it most} of them
$\Isom(M) \cong \{ e \}$ (our proof gives a partial support to believe that
this is the case).

\subsection{}
Another natural question is if a result like Theorem~1.1 will hold if we replace
$\Hy^n$ by other space $X$, say $X$ is $H/K$ where $H$ is a simple Lie group
and $K$ is a maximal compact subgroup of $H$.

In case $\R$-rank$(H) \ge 2$ one cannot expect this to be true. Indeed,
by Margulis' Theorem [Mr, Theorem~1, p.~2] every lattice $\Gamma$ in $H$ is
arithmetic, moreover, by Serre's conjecture~(cf. [PR, Section~9.5])
we expect $\Gamma$
to have the congruence subgroup property (in fact, Serre's conjecture has by
now been established for most of the cases). This gives a strong restriction
on the finite groups that can appear as quotients of finite index
subgroups of such $\Gamma$'s. For example, their Lie type composition factors
should have a bounded Lie rank depending only on $H$ and not on~$\Gamma$.

An analogue of Theorem~1.1 might hold for the complex hyperbolic
spaces $\Hy\C^n = \SU(n,1)/K$. Unfortunately, very little is currently known
here. In~\cite{Li}, Livne produced an example of a cocompact lattice in $\SU(2,1)$
which is mapped onto a non-abelian free group. This implies that for every
finite group $G$ there exists a compact manifold $M$ covered by $\Hy\C^2$, with
$G\subset\Isom(M)$. For $n>2$ we can not prove even this weak result.

\end{document}